\documentclass[a4paper]{article}
\usepackage[cp1251]{inputenc}
\usepackage[russian,english]{babel}
\usepackage{graphics}
\usepackage{amsmath,amsfonts,amssymb}

\newtheorem{theorem}{Theorem}[section]
\newtheorem{lemma}{Lemma}[section]
\newtheorem{corollary}{Corollary}[section]

\DeclareMathOperator{\Si}{Si} 
\def\myref#1{(\ref{#1})}
\def\R{{\mathbb R}}
\def\Z{{\mathbb Z}}

\DeclareMathOperator{\sign}{sign}

\makeatletter
\renewcommand{\@biblabel}[1]{#1.\hfill}
\@addtoreset{equation}{section} \@addtoreset{equation}{section*} \makeatother

\title{Complexity of the minimum-time damping of a physical pendulum}
\author{Alexander Ovseevich 
}
\date{}

\begin{document}
\maketitle

\begin{abstract}\noindent We study the  minimum-time damping of a physical pendulum by means of a
bounded control. In the similar problem for a linear oscillator each optimal trajectory possesses a finite number of
control switchings from  the maximal to the minimal value. If one considers simultaneously all optimal trajectories with
any initial state, the number of switchings can be arbitrary large. We show that for the nonlinear pendulum there is a
uniform bound for the switching number for all optimal trajectories. We find asymptotics for this bound as the control
amplitude goes to zero.
\end{abstract}

\medskip
\noindent{\bf Keywords:} bang-bang control; Sturm theory; Poincar\'e map

\medskip
\noindent{\bf MSC:} 49J15, 34C10


\section{Introduction}
The problem of minimum-time damping of a  pendulum is a classical issue of control theory. In the linear case, described
by the equation $\ddot x+x=u,\,|u|\leq1$, its solution is stated in \cite{pont}. The optimal control is of bang-bang
type, i.e. it takes values  $u=\pm1$, and the switching curve  separating the domain  $u=-1$ of the phase plane from the
domain $u=+1$ consists of unit semicircles centered at points of the form $(2k+1,0)$, where $k$ is an integer. The real
physical pendulum controlled by a torque in the joint is governed by the equation  $\ddot x+\sin x=\varepsilon
u,\,|u|\leq1$, where $x$ is the vertical angle, and $\varepsilon$ is the maximal amplitude of the control torque. The
parameter $\varepsilon$ is arbitrary: it can be large, small, of order 1. We are interested most in the case of a small
$\varepsilon$. The maximum principle says that the optimal control has the form $u=\sign\psi$, where the ``adjoint''
variable satisfies the equation $\ddot \psi+(\cos x)\psi=0$. Thus, the control is still of the bang-bang type, but the
time instants of switchings are the roots of a rather nontrivial function, a solution of the general
Sturm--Liouville/Schr\"odinger equation. The complexity of a control is characterized mainly by the switching number. In
the linear case this number for a trajectory connecting the initial point $(x,\dot x)$ with $(0,0)$ is
$\frac{T}{\pi}+O(1)$, where  $T$ is the duration of the motion. In turn, $T={\pi}\sqrt{\frac{E}{2}}+O(1)$, where
$E=\frac12{\dot x}^2+\frac12{x}^2$ is the energy of the system (cf. section \ref{damping}). Thus, each trajectory
possesses a finite number of switchings, but if the initial energy is large this number is $\sqrt{\frac{E}{2}}+O(1)$ and
is also large.

In the nonlinear case the switching number behaves quite differently. The best result, known to the author, is due to
Reshmin \cite{resh}. It says that if the parameter $\varepsilon$ is large enough, all optimal trajectories possess no
more than a single switching. Other interesting results can be found in \cite{bel,flag,PG}.

 We show that for any $\varepsilon$ the switching numbers for all optimal trajectories possess a
common bound. In other words, the following holds:

\begin{theorem}\label{fin} Suppose $N_\varepsilon(x,\dot x)$ is the number of zeroes of the adjoint variable
$\psi$  along an optimal trajectory connecting $(x,\dot x)$ with (0,0). Then the quantity $N_\varepsilon=\sup
N_\varepsilon(x,\dot x)$, where $\sup$ is taken over the entire phase space, is finite.
\end{theorem}
Another result states  upper and lower bounds for $N_\varepsilon$ which are  sharp with respect to the order of
magnitude.

\begin{theorem}\label{th2} There exist positive constants  $c_1,\,c_2$, 
such that
$$\frac{c_1}{\varepsilon}\leq N_\varepsilon\leq \frac{c_2}{\varepsilon}$$  for
$\varepsilon$  small enough.
\end{theorem}
Our main result is a promotion of inequalities of theorem \ref{th2} to an asymptotic equality:
\begin{theorem}\label{main0}There is the asymptotic equivalence
\begin{equation}\label{conj0}
 N_{\varepsilon}\sim \frac{D}\varepsilon \mbox{as }\varepsilon\to0,\mbox{ where }
D=\frac12\Si(\pi)=\int\limits_0^\pi\frac{\sin x}{2x}\,dx=0.925968526\dots
\end{equation}
\end{theorem}
The theorem can be regarded as an asymptotic formula $\varepsilon_n\sim D/n$ for bifurcation values of the parameter
$\varepsilon$. Here, the bifurcation is the  increment by 1 of the maximal number of the control switchings. The paper as
a whole grew, like the  ``Feigenbaum universality'' for the period-doubling bifurcation of a one-dimensional map
\cite{F}, out of contemplation of numerical data. This time the data were gathered  by S. Reshmin, who has computed the
complete phase portrait of the  minimum-time feedback control for many values of the parameter $\varepsilon$. In
particular, he has found the first 17 bifurcation values of $\varepsilon$ that conform to the theoretical value
\myref{conj0} of $D$ with  5-digit precision.

\medskip
\noindent The paper is based on a lemma saying that in the large speed area an optimal trajectory possesses no more than
a single switching. We use heavily the Sturm theory of root location for solutions of a Sturm--Liouville equation. It
allows us to relate $N_\varepsilon$ to the optimal time of motion from points with energy of order 1 to the point (0,0).
The lower bound given in Theorem \ref{th2} is based on energy considerations, which allow us to estimate this time. The
upper bound is more complicated and follows from a computation of the elapsed time in a motion under a quasioptimal
control. The asymptotic equivalence \myref{conj0} stems from the idea of the Poincar\'e map control, coupled with a
special nonlinear Sturm-like theorem.

\section{Problem Statement and Main Results}
We start with a more precise statement of the problem. The phase space is the tangent bundle $T(S^1)=S^1\times\R$ of the
circle, with coordinates $(x,y)$, where $x\in S^1=\R/2\pi\Z$, $y\in\R$. Physically speaking  $x$ is the vertical
deviation angle, so that the height of the pendulum over the horizontal plane is $1-\cos x$, and $y$ is the angular
velocity. The control system takes the form
\begin{equation}
\left\{
  \begin{array}{lll}\label{system}
\dot x&=&y,\\
\dot y&=&-\sin x+\varepsilon u,\, |u|\leq1
\end{array}\right.
\end{equation}
Under control $u=u(t)$ the system \myref{system} is a Hamiltonian one, and the corresponding Hamiltonian function is
$h(x,y)=\frac12{y}^2+(1-\cos x)+\varepsilon ux$. Here, of course, the phase space is not the cylinder $T(S^1)$, but the
covering plane with coordinates $(x,y)\in\R^2$.

We are interested in the minimum-time damping, i.e., the fastest motion from a given point $(x,y)\in T(S^1)$ to the
stable equilibrium  (lower) point $(0,0)$. An optimal control always exists: It suffices to show that every initial state
can be driven to the origin by an admissible control (i.e., that the system is controllable). Indeed, it is well known
(see \cite{filippov} Thm. 1) that if the set of admissible paths is not empty, the control system satisfies some
regularity and growth conditions, and the set of admissible velocities at any point of the phase space is convex, then
the motion time attains its infimum over the set  of admissible paths. It remains to establish controllability, for the
other conditions are trivially met. To the best of my knowledge, there is no general theorem implying this result. We
present an ``ad hoc'' proof in Section \ref{controllability}, where besides arguments of a wider applicability, a
specific ``dry-friction'' control is utilized.

According to the Pontryagin maximum principle this problem is associated with adjoint variables $(\phi,\psi)$ and the
Pontryagin function (Hamiltonian) $$H=y\phi+(-\sin x+\varepsilon u)\psi-1$$ so that to an optimal control $u$ that
maximizes the Pontryagin function there corresponds an optimal motion governed by  the corresponding canonical system,
and $H\equiv0$ along the optimal trajectory. In other words, besides the system \myref{system} the following relations
hold:
\begin{eqnarray}
&&u=\sign\psi,\\
&&\label{conj_system1}\dot\phi=(\cos x)\psi,\\
&&\label{conj_system2}\dot\psi=-\phi\\
&&\label{H} y\phi+(-\sin x+\varepsilon u)\psi-1\equiv0.
\end{eqnarray}
It follows immediately from \myref{conj_system1}, \myref{conj_system2} that a singular control or chattering are
impossible. The zeroes of the function $\psi$ cannot accumulate: in the limit point the vector $(\phi,\psi)$ vanishes,
which is incompatible with \myref{H}. Much more precise information is provided by the Sturm theory. When applied to the
equation $\ddot \psi+(\cos x)\psi=0$ it says that the distance between zeroes of  $\psi$ is no less than the distance
between zeroes of a solution to $\ddot \Psi+\Psi=0$, i.e. no less than $\pi=3.141\dots$ This follows just from the
inequality $\cos x\leq1.$\footnote{The Sturm theorem states \cite{codd_levinson} that if $f_i,\, i=1,2$ are solutions to
Sturm--Liouville equations $\ddot f_i+p_if_i=0$, and ``potentials'' $p_i$ are related by the inequality $p_1\leq p_2$,
then there is a root of the function $f_2$ between any pair of roots of $f_1$.} More generally one can state the
following corollary of \myref{conj_system1}, \myref{conj_system2}, and the Sturm theory:
\begin{lemma}\label{estim} In an optimal arc of duration $T$ no more than
$\frac{T}{\pi}+1$ switchings of control are
possible.
\end{lemma}
Note that the duration of an optimal motion can be arbitrary large, if the initial energy is large enough. Thus, the
Lemma does not immediately imply the finiteness Theorem \ref{fin}.

\subsection{Basic Lemma} We begin with a lemma which implies the  
absence of switchings
at high energy states.
\begin{lemma}[Basic Lemma]\label{switch} Suppose $t_1,\,t_2$ are adjacent zeroes of the adjoint
variable $\psi=\psi(t)$. Then, the velocities $y(t_1),\,y(t_2)$ have opposite directions.
\end{lemma}
\begin{corollary}\label{switch2}
Under conditions of the Lemma there is a time  $t$ between $t_1,\,t_2$ such that $y(t)=0.$
\end{corollary}
Notice that the statement of the Corollary is a Sturm-like theorem. In section \ref{sturm} we will prove a strengthening
of Corollary \ref{switch2}, where the uniqueness of  $t$ is asserted.

\medskip\noindent  The Corollary is obvious. To prove  Lemma \ref{switch} consider values of
$\phi(t_1),\,\phi(t_2)$. We assert that they have opposite signs. Indeed, $t_1$ and $t_2$ are adjacent simple zeroes of
$\psi,$ and the derivatives $\dot\psi(t_1),\,\dot\psi(t_2)$ should have opposite signs. Now, our statement follows from
\myref{conj_system2}. To complete the proof of the lemma we use equation \myref{H}, which implies that $y\phi=1$ at the
switching point, where $\psi=0.$ $\blacktriangleright$

The energy  \begin{equation}\label{energy0}E=\frac12{y}^2+(1-\cos x)\end{equation} of the pendulum cannot be large at a
point, where $y=0$; more precisely, at that point $E=|E|\leq2.$ Therefore, we get a corollary of the Basic Lemma to the
effect that even before the second switching the optimal motion takes place in the bounded energy area.

\medskip\noindent{\bf Remark.} Lemma \ref{switch} is not new \cite{lee_marcus}, but, to the best of my knowledge,
it has never been used to estimate the switching number for damping of a pendulum.

\subsection{Bounds for the damping time}
Let  $K$ be a compact in the phase space $T(S^1)$, and $T_{\varepsilon}=T_{\varepsilon}(K)$ be the maximum of damping
times over all initial conditions $(x,y)\in K$. Assume that $K$ is not the singleton (0,0). The  estimate for the time
$T_\varepsilon$ obtained in the next result provides the  ground for the finiteness Theorems \ref{fin} and \ref{th2}:

\begin{theorem}\label{th3} There exist positive constants  $C_i=C_i(K),\,i=1,2$ such that
$$\frac{C_1}{\varepsilon}\leq T_\varepsilon\leq \frac{C_2}{\varepsilon}$$ as the
{\em(}positive{\em)}  $\varepsilon$ is small enough.
\end{theorem}
Taking this  for granted, one can immediately prove Theorems  \ref{fin} and \ref{th2}. Indeed, we obtain from Lemmas
\ref{estim} and \ref{switch} that $N_{\varepsilon}\leq1+(\frac{T_{\varepsilon}(K)}{\pi}+1),$ where $K=\{E\leq2\},$ and
the energy $E=\frac12{y}^2+(1-\cos x).$ It is obvious that the damping time  $T_{\varepsilon}(K)$ is a monotone
decreasing function of $\varepsilon.$  Now, by taking $c_2=\frac{C_2(K)}{\pi}+1$ we immediately get Theorem \ref{fin} and
the upper estimate in Theorem \ref{th2}. To prove the lower estimate in Theorem \ref{th2} consider an optimal trajectory
connecting the oval ${\mathcal O}=\{E=1/2\}$ with the lower equilibrium state $(0,0)$, and contained in the set
$\widetilde{\mathcal O}=\{E\leq1/2\}.$ Such a trajectory should exist. Indeed, one can consider an optimal trajectory
running out of a point of the oval ${\mathcal O}$. From some instant on it is contained inside $\widetilde{\mathcal O}.$
At the instant mentioned the trajectory intersects the oval ${\mathcal O}$. We take the intersection point as the compact
$K$. Since $E\leq1/2$ along the trajectory, it follows that $1/2\leq\cos x$. Therefore, by applying the Sturm theory to
the equation $\ddot \psi+(\cos x)\psi=0$ we obtain that the number $N(\psi)$ of zeroes of the function $\psi$ on the
trajectory is no less than the number $N(\Psi)-1$, where $\Psi$ is a solution of $\ddot \Psi+(1/2)\Psi=0$. Of course
$N(\Psi)\geq\frac{T}{\sqrt{2}\pi}-1,$ where $T$ is the duration of the motion. According to Theorem \ref{th3} we obtain
that $T\geq{C_1}/{\varepsilon}$, and this proves the lower estimate for $N_{\varepsilon}$ in Theorem \ref{th2}, for it is
clear that
$N_{\varepsilon}\geq N(\psi).$ 

\subsection{Proof of Theorem \ref{th3}: lower estimate}
To prove the lower estimate it suffices to take a singleton for the compact  $K$. Take this point $p_0$ as the initial
one of an optimal trajectory $p(t)=(x,y)(t)$, and consider the energy $E$ of the running point on the trajectory as a
function of time. We have
\begin{equation}\label{energy_diff} \dot E=\varepsilon yu,\,|y|\leq\sqrt{2E},
\end{equation}
which implies that $\left|\frac{d}{dt}{\sqrt{E}}\right|=\left|\frac12\frac{\dot
E}{\sqrt{E}}\right|\leq\frac1{\sqrt{2}}\varepsilon $. Since the initial value of energy is $E(p_0)$, and the final one is
zero, we get a lower estimate for the elapsed time $T_\varepsilon\geq \frac{\sqrt{2E(p_0)}}{\varepsilon}$.
$\blacktriangleright$
\subsection{Proof of Theorem \ref{th3}: upper estimate}\label{upper}

This is the most complicated part of the paper, based on nontrivial estimates for the damping time for a given initial
state.

\noindent{\bf General strategy.} We divide the phase space into three parts: of high energy $\{E>2\}$, of low energy
$\{E<2\}$, and the standstill zone $S_{2\varepsilon}=\{|\sin x|<2\varepsilon,\,|y|<2\varepsilon\}$. For small
$\varepsilon$ the standstill zone consists of two connected components, the neighborhoods of the upper and lower
equilibrium point. To estimate the damping time we use a particular ``quasioptimal'' control which is given by the
formula
\begin{equation}\label{quasi}
   u=-\sign y
\end{equation} outside the  standstill zone. This control provides a kind of dry friction, and the  standstill
zone is the place, where the dry friction prevents any motion. The control reflects the idea of steepest local energy
descent. Note that on an interval of a constant velocity sign the controlled motion is governed by the Hamiltonian of the
form $\frac12{y}^2+(1-\cos x)\pm\varepsilon x.$ We will show that one can make it to the upper standstill zone from a
high energy state in time of order $O(1/\varepsilon)$, make it to the lower standstill zone from a low energy state in
time of the same order $O(1/\varepsilon)$, make it to the low energy zone from the upper standstill zone in time of order
$O(\log 1/\varepsilon)$, and, finally, make it to  the lower equilibrium point from the lower standstill zone in time of
order $O(1)$.

\medskip
\noindent{\bf Standstill zones.}

In order to understand the motion inside and in the vicinity of the standstill zone we use linearization of the control
system in a neighborhood of an equilibrium point. The linearization of nonlinear system \ref{system} at each equilibrium
point is a completely controllable linear system, if the control constraint $|u|\leq1$ is neglected. When the constraint
is imposed,  the system is just locally controllable. Therefore, it is possible in time of order $O(1)$ to pass along
feasible trajectories from any point of a disk centered in the  equilibrium of radius $c\varepsilon$ to any other point
of the disk. Here, $c>0$ is an absolute constant.  This follows in a formal way from the general result pertaining to an
arbitrary dimension:
\begin{theorem}\label{linearization} Suppose $\dot z=f(z)+\varepsilon g(z)u,\,|u|\leq1$ is a
controlled $C^2$-system in a neighborhood of $0\in\R^n,$ $f(0)=0$, and the corresponding linearized system $\dot
z=\frac{\partial f}{\partial z}(0)z+\varepsilon g(0)u$ is  completely controllable. Then, if the parameter $\varepsilon$
is small enough, the  set reachable from zero in time $\leq1$ along trajectories of the nonlinear system,   contains a
ball centered at zero of radius $c\varepsilon,$ where $c$ is a positive constant .
\end{theorem}
The theorem is well known, and follows, e.g, from Theorem 7 on p. 126 of \cite{sontag}. Unfortunately, the local
controllability does not guarantee that one can pass the standstill zone in time of order $O(1)$. Sometimes this is
impossible, but it is always possible to pass the standstill zone in time of order $O(\log 1/\varepsilon)$ by using
non-local maneuvers.

\medskip

\noindent  To this end we use the classical logarithmic bound for the oscillation period of the uncontrolled pendulum. It
has the following form. Let $p$ be a point of the phase space, denote by $\tau(p)$ the time required for the next hit of
the point $p$ in the  uncontrolled motion of the pendulum. Then, if the energy $E(p)=2+h,$ then
$\tau(p)=O(\log\left|\frac1h\right|)$ as $h\to0.$ Note that if $h=0$ the pendulum might stay forever in the upper
equilibrium state, so that $\tau=\infty$. Analytically, the estimate has the form
\begin{equation}\label{period}  \int_{0}^{2\pi}|\cos
s+1+h|^{-1/2}ds=O(\log |h|^{-1}) \mbox{ as } h\to0.
\end{equation}
A proof is presented in section \ref{period-estimate}. The elliptic integral in the left-hand side of \myref{period} is
deeply studied  from XVIII century on.

The manner of passage of the standstill zones is different in the upper and lower parts. The situation in the lower part
is simpler: The corresponding linearized system is globally controllable in spite of  the control bound $|u|\leq1$. This
follows, e.g., from the Brammer theorem \cite{brammer}. In fact, a stronger statement is proved in section \ref{damping}.
Therefore, the reachable set of the linearized system, starting from zero, in a sufficiently large time of order $O(1)$
contains any  disk centered at zero of radius  $C\varepsilon$, where $C>0$ is any given constant. The same holds true for
the nonlinear system, because the corresponding reachable sets differ wrt the Hausdorff metric in less than
$O(\varepsilon^2).$ Indeed, for any admissible control $|u(t)|\leq1$, where $t=O(1)$, all trajectories of both the linear
and the nonlinear system, starting at zero, stay within a disk of radius  $C\varepsilon$. The nonlinear part of the
right-hand side of differential equation is of order $O(\varepsilon^2)$ inside the disk. This implies that the difference
between solutions of the linear and the nonlinear equation is of the same order $O(\varepsilon^2)$. From these
considerations we conclude that it is possible to reach the lower equilibrium point from any point of the lower
standstill zone in time of order $O(1)$.

In order to get from a point  $p$ of the upper standstill zone to the low-energy zone one can do as follows: If the point
$p$ is at the distance less than $c\varepsilon$ from the  upper equilibrium point, we can move it in time  $O(1)$ to any
point at the distance  exactly $c\varepsilon$ from the upper equilibrium. This can be done by virtue of Theorem
\ref{linearization}. If $p$ is at the distance more than $c\varepsilon$ from the  upper equilibrium it stays intact.
Thereafter we switch the control off, and wait for the time  $t_1$, when the $x$-coordinate of the point $p(t_1)$ becomes
zero. Then we apply the control \myref{quasi} up to the time $t_2$, when the velocity $y$ of the point $p'=p(t_2)$
becomes zero. The energy decrease $E(p)-E(p')$ has the value $(\pi+o(1))\varepsilon$. Therefore, if
$E(p)=2+O(\varepsilon^2),$ then, we obtain $E(p')<2-\frac12\pi\varepsilon,$ provided that the parameter $\varepsilon$ is
sufficiently small. This means that the point  $p'$ is within the low-energy zone. In view of the estimate \myref{period}
for the period of oscillations the maneuver takes time of order $O(\log 1/\varepsilon)$.

\medskip

\noindent In what follows we do not consider the motion within and in the vicinity of the standstill zone. To estimate
the duration of motion within high and low energy zones we use the Poincar\'e section technique coupled with the
logarithmic bound for the oscillation period of the uncontrolled pendulum.

\medskip

\noindent{\bf High energies.} Consider first the high energy case. We estimate the damping time for the pendulum with
initial position $p$ of energy $E(p)=2+h,\,h>\varepsilon^2$. Throughout time $\leq\tau(p)=O(\log\frac{1}{\varepsilon})$
we use zero control until the coordinate $x$ becomes equal to $\pi$. Thereafter we apply control \myref{quasi}, and
consider the corresponding Poincar\'e map, i.e. we consider controlled trajectory up to the  time of hitting $x=\pi$. If
the condition
\begin{equation}\label{cond}
  h>2\pi\varepsilon
\end{equation} holds, the time will come. Otherwise, in time of order $O(\log\frac{1}{\varepsilon})$
we will arrive along the controlled trajectory at a zero speed position, which belongs either to the standstill zone or
the low energy zone  $E<2-\varepsilon^2.$ Suppose that condition \myref{cond} holds. Then, the  Poincar\'e map is
defined, and maps the point $(\pi,y)$ to $(\pi,y'),$ where
\begin{equation}\label{poincare}
    \frac{y'^2}{2}=\frac{y^2}{2}-2\pi\varepsilon.
\end{equation}
The duration of the controlled motion $\tau_\varepsilon$ is bounded from above by the duration $\tau$ of uncontrolled
oscillation from the new point $(\pi,y')$
\begin{equation}\label{transition}
   \tau_\varepsilon=\int\limits_0^{2\pi}\frac{dx}{\sqrt{2(h+1+\cos x-\varepsilon x)}}\leq
\int\limits_0^{2\pi}\frac{dx}{\sqrt{2(h+1-2\pi\varepsilon+\cos x)}}=\tau(\pi,y').
\end{equation}
In view of \myref{period} we get the upper bound \begin{equation}\label{tau-estimate}
    \tau\leq C\max\left\{\log\left|\frac1{h'}\right|,1\right\},
\end{equation}
where $h'=\frac{y'^2}{2}=h-2\pi\varepsilon$. We will repeat this energy lowering procedure
$h_{n+1}=h_{n}-2\pi\varepsilon$ while the condition \myref{cond} holds for $h=h_{n}$.  By virtue of \myref{poincare} the
number of steps is finite ($\leq\frac{h}{2\pi\varepsilon}$, where $h$ is the initial value of the reduced energy), and
the the total motion time is bounded from above, due to \myref{tau-estimate}, as follows
\begin{equation}\label{time}
T\leq C\left(\sum_{h_n\leq1}\log\frac1{h_n}+\sum_{h_n\geq1}1\right)\leq\frac{C}{2\pi\varepsilon} \int\limits_{0}^{1}
\log\frac1{x}\,dx+\frac{Ch}{2\pi\varepsilon}.
\end{equation}
Here, $h_n, \,n\geq1$ is the monotone decreasing sequence of the reduced energy values $E(p_n)-2$ at points obtained by
iteration of the Poincar\'e map \myref{poincare}. The right-hand side inequality follows from the monotonicity of the
logarithm:\begin{equation}\label{monotonicity}
2\pi\varepsilon\log\frac1{h_{n}}=(h_{n}-h_{n+1})\log\frac1{h_{n}}\leq\int_{h_{n+1}}^{h_n}\log\frac1{x}\,dx.
\end{equation}
Since the integral $\int_{0}^{1} \log\frac1{x}\,dx$ happily converges, the total duration of the iteration process is
$T=O(1/\varepsilon).$ Besides, some time is required to bring the point in a standard position at the initial and
terminal instant. The total motion time in the high energy zone is  $T+O(\log 1/\varepsilon)=O(1/\varepsilon).$

\medskip
\noindent{\bf Low energies.} Similarly, one can work out the case of low energies $E(p)=2-h,\,h>\varepsilon^2$ by using
the Poincar\'e map associated with the Poincar\'e section $\{y=0\}$ instead of $\{x=\pi\}.$ We confine ourselves without
loss of generality to the case, where the angle coordinate of the initial point $p$ is contained in the interval
$(-\pi,0),$ and the initial speed is zero: $p=(-x,0),\,x\in (0,\pi)$. One can bring any point in this position by using
uncontrolled motion throughout time $\leq\tau(p)=O(\log{1}/{\varepsilon})$ which does not affect  our final estimate
$T=O(1/\varepsilon)$ of the hitting time of the lower standstill zone. The control  \myref{quasi} forces the initial
point $p=(-x,0),\,x\in (0,\pi)$ to move along a trajectory of the canonical system with the Hamiltonian
$\frac12{y}^2+(1-\cos x)+\varepsilon x$ up to reaching a point $p'=(x',0),\,x'\in (0,\pi)$ of zero speed. There arises
the  Poincar\'e map $ p\mapsto p',$ where the initial and final energies are related by
\begin{equation}\label{energy}
  E(p)-E(p')=\varepsilon(x+x').
\end{equation}
In particular, for a constant  $c>0$ we have a lower bound for the energy loss
\begin{equation}\label{energ-ineq}
   E(p)-E(p')\geq c\varepsilon\sqrt{E(p)}.
\end{equation}
Totally similar relations hold for initial positions of the form $p=(x,0),\,x\in (0,\pi)$, which are obtained by symmetry
wrt the vertical axis. The transition time from $(-x,0)$ to $(x',0)$ is
\begin{equation}\label{transition3}
\tau_\varepsilon=\int\limits_{-x}^{x'}\frac{ds}{\sqrt{2(\cos s-\varepsilon s-\cos x)}}\leq
\int\limits_{-x}^{x'}\frac{ds}{\sqrt{2(\cos s-\varepsilon (x+x')-\cos x)}}.
\end{equation}
According to  \myref{period} the integral in the right-hand side is $O(\log 1/h'),$ where $h'=2-E(p')=1+\cos x',$ so that
\begin{equation}\label{ineq}
   \tau_\varepsilon\leq C\max\left\{\log\left|\frac1{h'}\right|,1\right\}.
\end{equation}
Let  $h_n, \,n\geq1$ be the monotone increasing sequence of the reduced energy values $2-E(p_n)$  at points obtained by
iteration of the Poincar\'e map $p_{n}\mapsto p_{n+1}$. By virtue of \myref{energy} we will make it to the lower
standstill zone in a finite number of steps. Moreover, the total time of the motion is estimated from above, thanks to
\myref{ineq}, as
\begin{equation}\label{time2}
T\leq C\sum_{h_n\leq1/e}\log\frac1{h_n}+C\sum_{h_n\geq1/e}1=C\sigma_1+C\sigma_2.
%
\end{equation}
In view of \myref{energ-ineq} we have $h_{n+1}-h_{n}\geq c\varepsilon\sqrt{2-h_{n}}$. If $h_n\leq1/e,$ where
$e=2.718\dots$ is the base of natural logarithms, this means that $h_{n+1}-h_{n}\geq c'\varepsilon.$ Therefore, the sum
$\sigma_1$ in the right-hand side of \myref{time2} can be estimated by $\frac1\varepsilon$ times the Riemann sum
\begin{equation}\label{time3}
\sum_{{h_n\leq1/e}}(h_{n+1}-h_{n})\log{h_n^{-1}}
\end{equation}
for the integral of logarithm over the interval $[0,1/e]$. Similarly, the sum $\sigma_2$ can be estimated via
$\frac1\varepsilon$ times the Riemann sum
\begin{equation}\label{time4}
\sum_{h_n\geq1/e}(h_{n+1}-h_{n})(2-h_{n})^{-1/2}
\end{equation}
for the integral of the function $x\mapsto(2-x)^{-1/2}$  over the interval $[1/e,2].$ The Riemann sums \myref{time3},
\myref{time4} can be estimated, like in \myref{monotonicity} via the convergent integrals
$$\int_{0}^{1/e} \log\frac1x\,dx\mbox{ and }\int_{1/e}^2 (2-x)^{-1/2}\,dx,$$ because $\log\frac1x$ and
$(2-x)^{-1/2}$ are monotone functions. Finally we arrive at the desired estimate $T=O(1/\varepsilon)$ for the hitting
time of the lower standstill zone.

\section{Theorem on the number of switchings}\label{conjecture}
The inequalities of Theorems  \ref{th2}, \ref{th3}  suggest the following natural question: Do there exist limits
${\mathfrak N}=\lim\limits_{\varepsilon\to0}\varepsilon N_{\varepsilon}$ and ${\mathfrak
T}=\lim\limits_{\varepsilon\to0}\varepsilon T_{\varepsilon}(K)$?  If they do, then, what are their values? Moreover, in
the proof we used a particular control $u$, namely, the one given by \myref{quasi} outside the standstill zones. A
natural question arises: how close is this control to the optimal one? For instance, is it true that the switching number
$ N_{\varepsilon}(u)$ and the damping time $T_{\varepsilon}(K)(u)$ for all points of a compact $K$ under this control are
asymptotically equivalent to that of the optimal control? More precisely, one can consider the limits ${\mathfrak
N}(u)=\lim\limits_{\varepsilon\to0}\varepsilon N_{\varepsilon}(u)$ and ${\mathfrak
T}(u)=\lim\limits_{\varepsilon\to0}\varepsilon T_{\varepsilon}(K)(u)$, and, if they exist, one may ask: Do they coincide
with ${\mathfrak N}$ and ${\mathfrak T}?$ To state the results we need to introduce the following notations. If $E\geq2$
put
\begin{equation}\label{tau+}
  \tau_+(E)=\frac1{2\sqrt{2}\pi}\int\limits_2^E dE\int\limits_0^{2\pi}(E-1-\cos\phi)^{-1/2}d\phi,
\end{equation}
and if $E\leq2$ let $\Phi=\Phi(E)\in[0,\pi]$ be the solution of $1-\cos\Phi=E$, and
\begin{equation}\label{tau-}
\tau_-(E)=\frac1{\sqrt{2}}\int\limits_0^\Phi \frac{\sin x}{x}\,dx\int\limits_0^x\frac{d\phi}{\sqrt{\cos \phi-\cos x}}.
\end{equation}
Finally, put \begin{equation}\label{tau}
   \tau(E)=\left\{\begin{array}{l}
\tau_+(E)+\tau_-(2),\mbox{ if }E\geq2\\
\tau_-(E),\mbox{ if }E\leq2
   \end{array}\right.
\end{equation}
\begin{theorem}\label{limit}
There exists the limit \begin{equation}\label{n-limit} {\mathfrak
N}(u)=\lim\limits_{\varepsilon\to0}\varepsilon\, N_{\varepsilon}(u)=
\int\limits_0^\pi\frac{\sin x}{2x}\,dx.
\end{equation} Suppose $E=E(K)=\max\limits_{p\in K}E(p)$ is the maximal energy of points of a
compact $K$. There exists the limit\begin{equation}\label{t-lim} {\mathfrak
T}(u)=\lim\limits_{\varepsilon\to0}\varepsilon\, T_{\varepsilon}(K)(u)=\tau(E),
\end{equation}
where $\tau(E)$ is defined in \myref{tau}.
\end{theorem}
Our main result is that the limits  ${\mathfrak N}$ and ${\mathfrak T}$ do exist and coincide with
those found in Theorem \ref{limit}. 

\begin{theorem}\label{main}Suppose
$N_{\varepsilon}$ is the maximal number of switchings for all mi\-n\-i\-mum-time trajectories of the system
\myref{system} ending at the lower equilibrium position $(0,0).$ Then, there exists the limit
\begin{equation}\label{conj}
{\mathfrak N}=\lim\limits_{\varepsilon\to0}\varepsilon N_{\varepsilon}\mbox{, and }{\mathfrak N}
=
\int\limits_0^\pi\frac{\sin x}{2x}\,dx=0.925968526\dots
\end{equation}
\end{theorem}
\begin{theorem}\label{t-main}Suppose
$T_{\varepsilon}(p)$ is the minimal time required for steering a point $p$ of the phase space of system \myref{system} to
the lower equilibrium position $(0,0),$ and $E=E(p)$ is the energy of $p$. Then, there exists the limit
\begin{equation}\label{t-conj}
  {\mathfrak T}(p)=\lim\limits_{\varepsilon\to0}\varepsilon\,
  T_{\varepsilon}(p)\mbox{, and }{\mathfrak T}(p)=\tau(E),
\end{equation}
where $\tau(E)$ is defined in \myref{tau}.
\end{theorem}
The proof of Theorem \ref{limit} is obtained by a refinement of the arguments used in section \ref{upper}. For instance,
the integrand $\frac{\sin x}{2x}$ in \myref{n-limit} arises as follows. We linearize the Poincar\'e map \myref{energy}
$$\cos x_{n+1}-\cos x_{n}=\varepsilon(x_{n}+x_{n+1}),$$ and obtain $-\frac{\sin
x_{n}}{2x_{n}}(x_{n}-x_{n+1})=\varepsilon+o(\varepsilon).$ In other words, we arrive at the Euler method with step
$\varepsilon$ for solution of the differential equation $-\frac{\sin x}{2x}\frac{dx}{dt}=1$, while the expression in the
right-hand side of \myref{n-limit} coincides with the time of motion from $x=\pi$ to 0. Furthermore, $\varepsilon
N_{\varepsilon}(u)$ is the discrete approximation of this time.

In order to prove Theorems \ref{main}, \ref{t-main} we need a more general Poincar\'e map control technique. The
technique can be regarded as a version of the averaging method well-known in the oscillation control
\cite{akulenko,akulenko2}. Besides, we use a strengthening of Corollary \ref{switch2} which allows us to count the
switchings via the number of changes of the speed direction.

\section{Poincar\'e map control}\label{poincare_control} In the previous sections \ref{upper} and
\ref{conjecture} we studied mainly a particular quasioptimal control \myref{quasi}. Still, some arguments can be applied
to an arbitrary control $u$.

To fix ideas, consider again the controlled motion of the point $p$ in the low energy zone $E(p)=2-h,\,h>\varepsilon^2$.
If a time interval under consideration is small compared to $1/\varepsilon,$  the trajectory $p(t)$ is close to the
trajectory of the uncontrolled motion with the same initial point. The Poincar\'e map ${\mathcal P}(u): p\mapsto p',$
related to the Poincar\'e section $\Sigma_-=\{y=0\},$ is close to the Poincar\'e map ${\mathcal P}(0)$ for  the
uncontrolled motion. In order to take into account the arising deviation of order $O(\varepsilon)$ it is convenient to
invoke  equation \myref{energy_diff} for energy change. Suppose $t_{n}$ are the hitting instances for the section
$\Sigma_-,$ $p_{n}=p(t_n)=(x_n,0)$ is the sequence of points arising under iteration of the Poincar\'e map,
$E_{n}=E(p_{n})$ are the corresponding values of energy. We have,\begin{equation}\label{en_diff}
  E_{n+1}-E_{n}=\varepsilon\int_{t_{n}}^{t_{n+1}} yudt=\varepsilon\int_{t_{n}}^{t_{n+1}} udx(t).
\end{equation}
We fix the time instant $t_{n}$, the point $p_{n}=(x_n,0)$, and study the influence of the control chosen upon the
right-hand side of \myref{en_diff}. In this equation $y=y(u,t)$ depends on control, but the effect is small:
$y(u,t)=y(0,t)+O(\varepsilon \tau_n),$ where $\tau_n=t_{n+1}-t_{n}$ is the time interval between next hits of the section
$\Sigma_-$. We know that in the  low energy zone $\tau_n=O(\log1/\varepsilon),$ and this bound is sharp in the vicinity
of the standstill zone only; in the major part of trajectory  $\tau_n$ is just bounded. If the time $t_{n}$ is fixed the
values of $t_{n+1}(u)$ and $\tau_{n}(u)$, like that of $y$, depend on $u$ weakly. Put
$$\phi_n=\int\limits_{t_{n}}^{t_{n+1}(0)} |y(0,t)|dt=\int\limits_{t_{n}}^{t_{n+1}(0)} |dx(0,t)|.$$
This is a function of the initial position $\phi_n=\phi(x_n).$ An easy computation shows that $\phi(x)=2|x|.$ Indeed, the
sign of the speed of uncontrolled pendulum does not change between adjacent positions with zero speed. Therefore,
\begin{equation}\label{phi}
  \int\limits_{t_{n}}^{t_{n+1}(0)} |dx(0,t)|=\left|\int\limits_{x_{n}}^{x_{n+1}^0} dx\right|=\left|x_{n+1}^0-x_{n}\right|,
\end{equation} where ${x_{n+1}^0}=x(0,t_{n+1}(0))$ is the
$x$-coordinate of the uncontrolled pendulum. It is clear that the right-hand side of \myref{phi} is $2|x_n|,$ because
under absence of control ${x_{n+1}^0}=-x_{n}.$
Thus, the right-hand side of 
\myref{en_diff} takes the form $\varepsilon\phi(x_n)U_n+o(\varepsilon),$ where $U_n$ is arbitrary subject to
$|U_n|\leq1.$ In the upshot we obtain a one-dimensional discrete control system
\begin{equation}\label{en_diff2}
  E_{n+1}-E_{n}=\varepsilon\phi(x_n)U_n+o(\varepsilon),\,|U_n|\leq1.
\end{equation}If we pass to the variables
$X_n=|x_n|$ the system takes the form
\begin{equation}\label{en_diff3}
\cos  X_{n}-\cos X_{n+1}=2\varepsilon X_{n}U_n+o(\varepsilon),\,|U_n|\leq1,
\end{equation}
or, equivalently,
\begin{equation}\label{en_diff4}
\frac{\sin  X_{n}}{2X_{n}}(X_{n+1}-X_{n})=\varepsilon U_n+o(\varepsilon),\,|U_n|\leq1.
\end{equation}
The obtained discrete  system arises via approximation by the Euler broken lines with step $\varepsilon$ of the
continuous  control system
\begin{equation}\label{en_diff5}
\frac{\sin  X}{2X}\frac{d  X}{dt}= U,\,|U|\leq1,
\end{equation}
so that $X_n$ approaches $X(n\varepsilon)$.

A similar  control system arise in the high energy zone. It is convenient to  $\Sigma_+=\{x=\pi\}$ as the Poincar\'e
section in the zone, and consider the controlled sequence   $y_n$ of velocities of the pendulum at hitting times $t_n$
for the section $\Sigma_+$. Then, the analogue of \myref{en_diff2} looks like
\begin{equation}\label{poincare+}
    \frac{y_{n+1}^2}{2}-\frac{y_{n}^2}{2}=2\pi\varepsilon U_n+o(\varepsilon),\,|U_n|\leq1.
\end{equation}
Of couurse, like in the low energy case, the discrete system arises out approximation by the Euler broken lines with step
$\varepsilon$ of the continuous control system
\begin{equation}\label{en_diff6}
Y\frac{d  Y}{dt}=2\pi U,\,|U|\leq1,
\end{equation}
where  $y_n$ approaches $Y(n\varepsilon)$.

The use of control  \myref{quasi} corresponds to $U\equiv-1.$ The minimum-time damping problem corresponds to
minimization of the functional $\sum_n \tau_n.$ After normalization $\sum_n \tau_n\mapsto \varepsilon\sum_n \tau_n$ and
passage to the limit $\varepsilon\to0$, if the initial position belongs to the low energy zone, we get the problem of
steering the system \myref{en_diff5} to the point $X(T)=0$ coupled with minimization of the functional
\begin{equation}\label{J-}
J_-=\int_0^T \tau_-(X(t))dt\to\min,\, 
\tau_-(X)= \int_0^{X} (\cos\phi-\cos X)^{-1/2}d\phi.
\end{equation}If the initial position belongs to the
high energy zone, there arises an extra problem of bringing system \myref{en_diff6} to  $Y(T)=0$ coupled with
minimization of the functional\begin{equation}\label{J+}
J_+=\int_0^T \tau_+(Y(t))dt,\to\min,\, 
\tau_+(Y)= \int_0^{2\pi} (1+Y^2/2 -\cos\phi)^{-1/2}d\phi.
\end{equation}It is more or less clear that the control $U\equiv-1$,
corresponding to \myref{quasi}, is optimal in both cases.
\subsection{Convergence of the Euler
broken lines}\label{euler} We present some details on convergence of the Euler broken lines for equations
\myref{en_diff5}, \myref{en_diff6}. They are not totally standard, because these equations are implicit.

We extend the discrete sequence \myref{en_diff3} to the piecewise-linear function $X_{\varepsilon}$ by defining
$n=[t/\varepsilon]$, and $X_{\varepsilon}(t)=X_{n}+(t-\varepsilon n)X_{n+1}$ if $t\in[\varepsilon n,\varepsilon (n+1)].$
We extend the sequence of controls $U_n$ to the fuction $U_{\varepsilon},$ which is constant in the intervals
$[\varepsilon n,\varepsilon (n+1)].$ The domain of the functions $X_{\varepsilon},\,U_{\varepsilon}$ is not known in
advance. It follows from the upper estimate  in Theorem \ref{th3} for the duration  of the optimal motion that for the
optimal control this domain is bounded.

The functions $f_{\varepsilon}=\cos X_{\varepsilon}$ form an equicontinuous family, because it follows from
\ref{en_diff3} that $|f_{\varepsilon}(t)-f_{\varepsilon}(s)|\leq C|t-s|.$ Therefore, there exists a subsequence
$\varepsilon=\varepsilon_k\to0$such that the functions $f_{\varepsilon}$ converge uniformly on bounded intervals to a
function $f.$ Therefore, the same is true for the functions $X_{\varepsilon}\to X$. By taking a subsequence one can
assume that the functions $U_{\varepsilon}$ converge to a function $U,\, |U|\leq1$ weakly, i.e. $\int
g(\sigma)U_{\varepsilon}(\sigma)d\sigma\to \int g(\sigma)U(\sigma)d\sigma$ for any fixed integrable function $g\in L_1.$
It follows from equation \ref{en_diff3} that as $\varepsilon\to0$
\begin{equation}\label{st}
  \cos  X_{\varepsilon}(s)-\cos X_{\varepsilon}(t)=2\int_s^t
  X_{\varepsilon}(\sigma)U_{\varepsilon}(\sigma)d\sigma+o(1),
\end{equation}
and, therefore, in the limit
\begin{equation}\label{stu}
  \cos  X(s)-\cos X(t)=2\int_s^t X(\sigma)U(\sigma)d\sigma.
\end{equation}
The integral equation \myref{stu} is equivalent to the differential equation \myref{en_diff5}. Since the differential
equation \myref{en_diff5} has a unique solution with initial condition $X(0)=\pi,$ the Euler broken lines
$X_{\varepsilon}$ converge uniformly as $\varepsilon\to0$, if the controls $U_{\varepsilon}$ converge to $U$ weakly.

Similarly one can obtain the convergence of the Euler broken lines $Y_{\varepsilon}$, where
$Y_{\varepsilon}(n\varepsilon)=y_n,$ in the high energy zone.
\section{One more ``Sturm-like'' theorem}\label{sturm}
To prove the main result on asymptotics of $N_\varepsilon$ we need the following  `Sturm-like'' strengthening of
Corollary \ref{switch2}.

\begin{theorem}\label{switch3}  Suppose $\varepsilon$ is sufficiently small, $t_1,\,t_2$
are adjacent zeroes of the adjoint variable $\psi=\psi(t)$, and the optimal motion in the interval $[t_1,t_2]$ of time
does not hit the standstill zone. Then, there exists a single time instant $t$ between $t_1,\,t_2$ such that $y(t)=0,$ so
that the zeroes of $y$ and the adjoint variable $\psi$ are intermittent.
\end{theorem}
Indeed, otherwise there are at least three subsequent time instants $\tau_1<\tau_2<\tau_3$  in the time interval
$[t_1,t_2]$ such that $y(\tau_i)=0,$ and the optimal control $u$ is constant in the interval $[\tau_1,\tau_3].$ To fix
ideas, assume that  $u=1$. Then, in the interval $[\tau_1,\tau_3]$ the Hamiltonian $\frac12y^2+(1-\cos x)+\varepsilon x$
is constant. In particular, the points $x_i=x(\tau_i)$ are in the same level of potential energy
$U_\varepsilon(x)=(1-\cos x)+\varepsilon x.$ If $\varepsilon=0$ the triple has the form $(x,-x,x)$. The perturbed triple
$x_i,\,i=1,2,3$ has the form $(x,-x',x''),$ where $x',x''$ are close to $x$ and equation $-\cos x+\varepsilon x=-\cos
x''+\varepsilon x''$ holds. We show, like in the implicit function theorem, that the only solution  close to  $x$ is
$x''=x.$ Indeed, if $x''\neq x$ the equation can be rewritten as $\frac{\cos x-\cos x''}{x-x''}=\varepsilon.$ By the mean
value theorem we obtain $\sin z=-\varepsilon,$ where $z$ is contained in the interval $[x,x''].$ Under condition $|\sin
x_i|\geq2\varepsilon$, which reflects that the motion goes outside the standstill zone, the latter equality is
impossible. Thus, the perturbed triple $x_i,\,i=1,2,3$ has the form $(x,-x',x)$. This means, however, that the points
$p(\tau_1)=(x_1,0)$, and $p(\tau_3)=(x_3,0)$ on the optimal trajectory coincide, which is impossible, because the elapsed
time $\tau_3-\tau_1$ is positive.
\section{Coda: proof of Theorems \ref{main}, \ref{t-main}} The proof is a combination of results
obtained in sections \ref{poincare_control}, and \ref{sturm}. Indeed, as $\varepsilon\to0$ the initial minimum-time
problem reduces to the optimal control problem for systems \myref{en_diff5}, \myref{en_diff6}, and moreover, according to
Theorem \ref{switch3}, the total time of motion for these new control systems corresponds to the limit of $\varepsilon
N_\varepsilon.$ The optimal control in the pendulum damping problem corresponds to the optimal control $U=-1$ in systems
\myref{en_diff5}, \myref{en_diff6}. For this control the computation of the time of motion is obvious, and is already
made within the proof of Theorem \ref{limit}.

\section{Auxiliary results}\label{auxiliar}
In this section we collect proofs of several already used auxiliary results.

\subsection{Controllability of the physical pendulum}\label{controllability}

Here, we prove that system \myref{system} is controllable. It suffices to show that zero is reachable from any point, or what is equivalent, that any point is
reachable from zero.

We first prove by rather general arguments that the low-energy zone $\{E<2\}$ is reachable from zero. Indeed, within this
zone all trajectories of the free motion, where $u=0$, are compact curve, which are energy levels.  The motion along any
of these curves is periodic. This implies, in particular, that the reachable part $R$ of the low-energy zone is a union
of the above curves. By applying control $u=\pm1$ at a suitable point of a curve we can always increase or decrease the
energy. Therefore, the reachable levels of energy inside the low-energy zone fill the open interval $(0,2)$. Therefore,
the entire low-energy zone is reachable. Since the closure of $R$ is also reachable we conclude that the entire closure
$\mathcal{E}=\{E\leq2\}$ of the  low-energy zone is reachable.

Second, we prove that from any point $(x_0,y_0)$ of high-energy zone $\{E>2\}$ one can reach the set $\mathcal{E}$, which
completes our arguments. To do this we apply the ``dry-friction'' control of the form $u=-\delta\sign y$ where $\delta$
is much less than $\varepsilon$. The control is used within the domain $\{\sign y=\sign y_0\}$. As soon as we hit the
boundary $\{ y=0\}$ we are in $\mathcal{E}$. Under the control, the energy decreases according to $\dot E=-\delta|y|$.
Therefore, in finite time we reach the standstill zone $S_{2\delta}=\{|\sin x|<2\delta,\,|y|<2\delta\}$. The standstill
zone is a neighborhood of two equilibrium points. By linearizing our system at an equilibrium point we'll get a
completely controllable linear system, and this implies (see Theorem \ref{linearization} or \cite{sontag} Theorem 7 on p.
126) that we can reach from our initial point an equilibrium in a finite time. Since the equilibria is contained in
$\mathcal{E}$, the proof is complete.

\subsection{Damping time in the linear problem}\label{damping}

Here, we prove a formula from the Introduction:  $T={\pi}\sqrt{\frac{E}{2}}+O(1),$   for the damping time of the linear
oscillator. Here, $E$ is the initial energy. The oscillator is governed by system
\begin{equation}
\left\{
  \begin{array}{lll}\label{lin-system}
\dot x&=&y,\\
\dot y&=&-x+u,\, |u|\leq1.
\end{array}\right.
\end{equation}
We have to show that the  set  ${\mathcal D}_T$, reachable from zero in time $T$, contains the disk ${\mathcal
E}=\left\{\frac12{y}^2+\frac12{x}^2\leq E\right\}$ for a sufficiently large $T$ related to $E$ via
$T={\pi}\sqrt{\frac{E}{2}}+O(1)$. The support function of the set ${\mathcal D}_T$ has the explicit expression
$$H_T(\xi)=\int_0^T\left|\xi_1\sin t+\xi_2\cos t\right|dt,$$ which is a particular case of the general
integral formula
$$\int_0^T\left|B^*e^{A^*t}\xi\right|dt$$ for support function of the reachable set of a linear
system $$\dot z=Az+Bu,\,|u|\leq1.$$  We show that
\begin{equation}\label{ineq-supp} H_T(\xi)=\frac2\pi|\xi|(T+O(1)),
\end{equation} where $|\xi|=(\xi_1^2+\xi_2^2)^{1/2}$ is the Euclidean norm of the vector $\xi.$
Since the support function for the disk  ${\mathcal E}$ equals $H_{\mathcal E}(\xi)=\sqrt{2E}|\xi|$, then, perhaps after
an increase of $T$ by a bounded value, we will get the required inequality $H_T\geq~H_{\mathcal E}.$ As it is well-known
the inequality between support functions implies inclusion of the corresponding closed convex sets. It remains to prove
\myref{ineq-supp}. We assume $|\xi|=1$ for simplicity. Consider the difference
$$\phi(T)=H_T(\xi)-\frac{2T}\pi=\int_0^T\left(f(t)-\frac2\pi\right)dt.$$
Here, $f$ is the  $\pi$-periodic function $|\xi_1\sin t+\xi_2\cos t|=|\cos(t-\alpha)|$,  where $\tg\alpha=\xi_1/\xi_2.$
The integral $I=\int\limits_{0}^{\pi}fdt$ over period equals $\int\limits_{0}^{\pi}|\cos t|dt=2$. Therefore, the function
$\phi$ is a $\pi$-periodic one. In particular, it is bounded, which proves \myref{ineq-supp}.

A minor extra effort allows us to obtain a precise estimate for $\Phi=\max\limits_{T,\xi}|\phi(T)|.$ Put
$\phi_0(T)=\int_0^T|\cos t|dt-\frac{2T}\pi$ corresponding to the value $\alpha=0$ of the shift parameter, and define
$\Phi_0=\max\limits_{T}\phi_0(T).$ Then, one can see that $\Phi=2\Phi_0,$ and the maximum of $\phi_0(T)$ is attained in
the interval $(0,\frac\pi2)$ when $\cos T=\frac{2}{\pi}.$ Therefore, the value of the sharp bound $2\Phi_0$ of the error
term in \myref{ineq-supp} can be made explicit:
$\Phi_0=\sqrt{1-\frac4{\pi^2}}-\frac{2}{\pi}\arccos\frac{2}{\pi}=0.2105\dots$.

\subsection{Estimate for the period of oscillation of a pendulum}\label{period-estimate}
{
Let us prove \myref{period}. It is clear that the only essential part of the integration interval is a small neighborhood
$(\pi-\delta,\pi+\delta)$ of $s=\pi,$ because the integral over the complement in $[0,2\pi]$ is bounded. We use in the
neighborhood the expansion $\cos(\pi+s)=-1+\frac{s^2}2+O(s^4).$  If $\delta$ is small enough, the term $O(s^4)$ in the
last equality does not exceed $\frac{s^2}4$ in absolute value, and, therefore,
$|\cos(\pi+s)+1+h|\geq\min\left\{\left|\frac{s^2}4+h\right|,\left|\frac{3s^2}4+h\right|\right\},$ which implies
$|\cos(\pi+s)+1+h|^{-1/2}\leq \left|\frac{s^2}4+h\right|^{-1/2}+\left|\frac{3s^2}4+h\right|^{-1/2}.$ It remains to obtain
the estimate of the form  $\int_{-1}^{1}|As^2+h|^{-1/2}ds=O(\log |h|^{-1}),$ where $A$ is an arbitrary positive constant.
The substitution $h\mapsto h/A$ reduces the task to the case $A=1.$ We write $h$ in the form $h=\pm B^2.$ The change of
variables $\sigma=s/B$ reduces  \myref{period} to the bound $\int_{-1/B}^{1/B}|\sigma^2\pm 1|^{-1/2}d\sigma=O(\log
|B|^{-1}),$ which is true, because the possible singularity  $\sigma=\pm1$ of the integrand is integrable, and at
infinity $|\sigma^2\pm 1|^{-1/2}\sim|\sigma|^{-1}$.}

\section{Conclusions}
The above results by no means give a clear picture of the complexity of the minimum-time problem for general control
systems. Which control systems possess a finite number of switchings? Where is the divide  between the minimum time
problem for a pendulum, and the Fuller problem \cite{borisov,zb}, where chattering is a stable phenomenon?

It is  clear that our results can be extended to the problem of bringing the pendulum to any given state, not just the
lower equilibrium position, and it seems  possible to extend the above finiteness theorems to systems governed by general
nonlinear second order differential equations $ \ddot x+f(x) =\varepsilon u,\,|u|\leq1.$ A further reasonable step would
be the study of finiteness phenomena for general 2-dimensional systems in the spirit of \cite{BP,sussmann}. It is
absolutely unclear, however, what's going on with  multidimensional systems, e.g. for two interrelated nonlinear
pendulums. One can claim a general conjecture to the effect that finiteness of the number of switchings holds for a
generic controlled Hamiltonian system defined on the cotangent bundle of a compact manifold. This is a rather bold step.
The statement itself requires a clarification, and many aspects of our methods, intimately related to peculiarities of
the pendulum, can hardly be extended to many dimensions. It seems that the most enigmatic is a multidimensional
counterpart of Corollary \ref{switch2} related to the phenomenon of almost complete absence of switchings outside a
compact part of the phase space.
\section*{Acknowledgments} This work was supported by  Russian Foundation of Basic Research, grant
11-08-00435.


\begin{thebibliography}{9}
\bibitem{pont} {\em Pontryagin L.S., Boltyansky V.G., Gamkrelidze  R.V., Mischenko E.F.} Mathematical theory
of optimal processes. Moscow: Nauka (1983)

\bibitem{resh} {\em Reshmin S.A.}//
Journal of Applied Mathematics and Mechanics 73(4) (2009) 403–410.

\bibitem{bel} {\em Beletsky V.V.}// Space Studies.  v. 9, no. 3. p. 366--375 (1971)

\bibitem{flag} {\em Garcia Almuzara J.L., Fl\"ugge-Lots I.} Minimum time control of a nonlinear system // J.
Differential Equations.  Vol. 4, no. 1, pp. 12--39 (1968)

\bibitem{lee_marcus} {\em Lee E. B., Marcus L.} On necessary and sufficient optimality conditions
in minimum-time problem for nonlinear second order systems// Proc. II Int. Congress IFAC. Basel,
(1963) 

\bibitem{codd_levinson} {\em Coddington E., Levinson N.} Theory of ordinary differential equations.
McGraw-Hill, New York (1955)

\bibitem{borisov}{\em Borisov V.F.}  Fuller's Phenomenon: Review, J. Math. Sciences, Vol. 100, No. 4,
(2000)

\bibitem{akulenko}{\em Akulenko L.D.} Asymptotic methods of optimal control, Moscow: Nauka (1987)

\bibitem{akulenko2}{\em Chernousko F.L., Akulenko L.D., Sokolov B.N.} Control of Oscillations. Moscow: Nauka,  383
pp. (1980)

\bibitem{brammer}{\em Brammer R.F.} Controllabibity of linear autonomous
systems with positive passive controllers.  SIAM J. on Control,  v. 10, No. 2, p. 339--353  (1972)

\bibitem{BP}{\em U. Boscain and B. Picolli,} Optimal Syntheses for Control Systems on 2-D Manifolds, Springer, Berlin, 2004 , Vol. 43.

\bibitem{sussmann}{\em  H. J. Sussmann.} Regular synthesis for time-optimal control of single-input real analytic
systems in the plane, SIAM J. Control and Optimization, Vol. 25, No. 5, September 1987

\bibitem{PG}{\em Paoletti, P. and Genesio, R.} Rate limited time optimal control of a planar pendulum, Systems
Control Lett., 60(2011), no. 4, 264–270.

\bibitem{sontag}{\em E. D. Sontag} Mathematical Control Theory (2nd Edition), Springer-Verlag, New York, 1998

\bibitem{zb}{\em Zelikin, M. I.; Borisov V. F.} Theory of chattering control. With applications to astronautics,
robotics, economics, and engineering. Systems \& Control: Foundations \& Applications. Birkha\"user
Boston, Inc., Boston, MA, 1994. 

\bibitem{F}{\em  Feigenbaum, M. J.} (1978). Quantitative Universality for a Class of Non-Linear Transformations. J. Stat. Phys. 19: 25–52.

\bibitem{filippov}{\em Filippov A. F.} On certain questions in the theory of
optimal control. Vestnik Moskov. Univ. Ser. Mat. Mech. Astr. 2 (1959) 25–32 (English translation: SIAM J. Control 1
(1962) 76–84.).
\end{thebibliography}
\end{document}